\documentclass[twoside]{ima}
\usepackage{amsmath}
\usepackage{amssymb}
\usepackage[all]{xypic}
\let\Cal\mathcal
\let\Bbb\mathbb
\let\phi\varphi

\newcommand{\x}{\times}
\renewcommand{\o}{\circ}
\newcommand{\al}{\alpha}
\newcommand{\be}{\beta}
\newcommand{\Ga}{\Gamma}
\newcommand{\La}{\Lambda}
\newcommand{\om}{\omega}
\newcommand{\Om}{\Omega}
\newcommand{\ph}{\phi}
\newcommand{\Ph}{\Phi}
\newcommand{\ps}{\psi}

\newcommand{\la}{\lambda}

\newcommand{\Si}{\Sigma}

\newcommand{\fg}{{\mathfrak g}}

\newcommand{\ad}{\operatorname{ad}}

\newcommand{\im}{\operatorname{im}}
\newcommand{\id}{\operatorname{id}}

\begin{document}
\markboth{ANDREAS \v CAP}{OVERDETERMINED SYSTEMS, CONFORMAL GEOMETRY
  AND BGG}

\title{OVERDETERMINED SYSTEMS,\\ CONFORMAL DIFFERENTIAL GEOMETRY,\\
  AND THE BGG COMPLEX} 

\author{ANDREAS \v CAP}\thanks{Fakult\"at f\"ur Mathematik,
  Universit\"at Wien, Nordbergstra\ss e~15, A--1090 Wien, Austria, and
  International Erwin Schr\"odinger Institute for Mathematical
  Physics, Boltzmanngasse 9, A-1090 Wien, Austria, e-mail:
  andreas.cap@esi.ac.at\\supported by project P15747--N05 of the
  ``Fonds zur F\"orderung der wissenschaftlichen Forschung'' (FWF) and
  by the Insitute for Mathematics and its Applications (IMA)}

\maketitle

\begin{abstract}
  This is an expanded version of a series of two lectures given at the
  IMA summer program ``Symmetries and Overdetermined Systems of
  Partial Differential Equations''. The main part of the article
  describes the Riemannian version of the prolongation procedure for
  certain overdetermined system obtained recently in joint work with
  T.P.~Branson, M.G.~Eastwood, and A.R.~Gover. First a simple special
  case is discussed, then the (Riemannian) procedure is described in
  general.
  
  The prolongation procedure was derived from a simplification of the
  construction of Bernstein--Gelfand--Gelfand (BGG) sequences of
  invariant differential operators for certain geometric structures.
  The version of this construction for conformal structures is
  described next. Finally, we discuss generalizations of both the
  prolongation procedure and the construction of invariant operators
  to other geometric structures.
\end{abstract}


{\AMSMOS 35N10, 53A30, 53A40, 53C15, 58J10, 58J70\endAMSMOS}

\section{Introduction}
The plan for this article is a follows. I'll start by describing a
simple example of the Riemannian version of the prolongation procedure
of \cite{BCEG}. Next, I will explain how the direct observations used
in this example can be replaced by tools from representation theory to
make the procedure work in general. The whole procedure is based on an
inclusion of the group $O(n)$ into $O(n+1,1)$.  Interpreting this
inclusion geometrically leads to a relation to conformal geometry,
that I will discuss next. Via the conformal Cartan connection, the
ideas used in the prolongation procedure lead to a construction of
conformally invariant differential operators from a twisted de--Rham
sequence. On manifolds which are locally conformally flat, this leads
to resolutions of certain locally constant sheaves, which are
equivalent to the (generalized) Bernstein--Gelfand--Gelfand (BGG)
resolutions from representation theory. In the end, I will outline
generalizations to other geometric structures.

It should be pointed out right at the beginning, that this
presentation basically reverses the historical development. The BGG
resolutions in representation theory were originally introduced in
\cite{BGG} and \cite{Lepowsky} in the 1970's. The constructions were
purely algebraic and combinatorial, based on the classification of
homomorphisms of Verma modules. It was known to the experts that there
is a relation to invariant differential operators on homogeneous
spaces, with conformally invariant operators on the sphere as a
special case. However it took some time until the relevance of ideas
and techniques from representation theory in conformal geometry was
more widely appreciated. An important step in this direction was the
work on the curved translation principle in \cite{Eastwood-Rice}. In
the sequel, there were some attempts to construct invariant
differential operators via a geometric version of the generalized BGG
resolutions for conformal and related structures, in particular in
\cite{Baston}.

This problem was completely solved in the general setting of parabolic
geometries in \cite{CSS:BGG}, and the construction was significantly
simplified in \cite{CD}. In these constructions, the operators occur
in patterns, and the first operators in each pattern form an
overdetermined system. For each of these systems, existence of
solutions is an interesting geometric condition. In \cite{BCEG} it was
shown that weakening the requirement on invariance (for example
forgetting conformal aspects and just thinking about Riemannian
metrics) the construction of a BGG sequence can be simplified.
Moreover, it can be used to rewrite the overdetermined system given by
the first operator(s) in the sequence as a first order closed system,
and this continues to work if one adds arbitrary lower order terms.

I am emphasizing these aspects because I hope that this will clarify
two points which would otherwise remain rather mysterious. On the one
hand, we will not start with some overdetermined system and try to
rewrite this in closed form. Rather than that, our starting point is
an auxiliary first order system of certain type which is rewritten
equivalently in two different ways, once as a higher order system and
once in closed form. Only in the end, it will follow from
representation theory, which systems are covered by the procedure. 

On the other hand, if one starts the procedure in a purely Riemannian
setting, there are some choices which seem unmotivated. These choices
are often dictated if one requires conformal invariance.

\section{An example of the prolongation procedure}\label{2} 
\setcounter{proposition}1
\setcounter{theorem}3
\subsection{The setup}\label{2.1}
The basics of Riemannian geometry are closely related to
representation theory of the orthogonal group $O(n)$. Any
representation of $O(n)$ gives rise to a natural vector bundle on
$n$--dimensional Riemannian manifolds and any $O(n)$--equivariant map
between two such representation induces a natural vector bundle map.
This can be proved formally using associated bundles to the
orthonormal frame bundle.

Informally, it suffices to know (at least for tensor bundles) that the
standard representation corresponds to the tangent or cotangent
bundle, and the correspondence is natural with respect to direct sums
and tensor products. A linear map between two representations of
$O(n)$ can be expressed in terms of a basis induced from an
orthonormal basis in the standard representation. Starting from a
local orthonormal frame of the (co)tangent bundle, one may locally use
the same formula in induced frames on any Riemannian manifold.
Equivariancy under the group $O(n)$ means that the result is
independent of the initial choice of a local orthonormal frame.  Hence
one obtains a global, well defined bundle map.

The basic strategy for our prolongation procedure is to embed $O(n)$
into a larger Lie group $G$, and then look how representations of $G$
behave when viewed as representations of the subgroup $O(n)$. In this
way, representation theory is used as a way to organize symmetries. A
well known inclusion of this type is $O(n)\hookrightarrow O(n+1)$,
which is related to viewing the sphere $S^n$ as a homogeneous
Riemannian manifold. We use a similar, but slightly more involved
inclusion.

Consider $\Bbb V:=\Bbb R^{n+2}$ with coordinates numbered from $0$ to
$n+1$ and the inner product defined by 
$$
\langle(x_0,\dots,x_{n_1}),(y_0,\dots,y_{n+1})\rangle:=
x_0y_{n+1}+x_{n+1}y_0+\sum_{i=1}^nx_iy_i. 
$$
For this choice of inner product, the basis vectors $e_1,\dots,e_n$
span a subspace $\Bbb V_1$ which is a standard Euclidean $\Bbb R^n$,
while the two additional coordinates are what physicists call light
cone coordinates, i.e.~the define a signature $(1,1)$ inner product on
$\Bbb R^2$. Hence the whole form has signature $(n+1,1)$ and we
consider its orthogonal group $G=O(\Bbb V)\cong O(n+1,1)$. There is an
evident inclusion $O(n)\hookrightarrow G$ given by letting $A\in O(n)$
act on $\Bbb V_1$ and leaving the orthocomplement of $\Bbb V_1$ fixed.

In terms of matrices, this inclusion maps $A\in O(n)$ to the block
diagonal matrix $\left(\begin{smallmatrix} 1 & 0 & 0\\ 0 & A & 0\\ 0 &
    0 & 1\end{smallmatrix}\right)$ with blocks of sizes $1$, $n$, and
$1$. The geometric meaning of this inclusion will be discussed later.

The representation of $A$ as a block matrix shows that, as a
representation of $O(n)$, $\Bbb V=\Bbb V_0\oplus\Bbb V_1\oplus V_2$,
where $\Bbb V_0$ and $\Bbb V_2$ are trivial representations spanned by
$e_{n+1}$ and $e_0$, respectively. We will often denote elements of
$\Bbb V$ by column vectors with three rows, with the bottom row
corresponding to $\Bbb V_0$.

If we think of $\Bbb V$ as representing a bundle, then differential
forms with values in that bundle correspond to the representations
$\La^k\Bbb R^n\otimes\Bbb V$ for $k=0,\dots,n$. Of course, for each
$k$, this representation decomposes as $\oplus_{i=0}^2(\La^k\Bbb
R^n\otimes\Bbb V_i)$, but for the middle component $\La^k\Bbb
R^n\otimes\Bbb V_1$, there is a finer decomposition. For example, if
$k=1$, then $\Bbb R^n\otimes\Bbb R^n$ decomposes as
$$
\Bbb R\oplus S^2_0\Bbb R^n\oplus\La^2\Bbb R^n
$$
into trace--part, tracefree symmetric part and skew part. We actually
need only $k=0,1,2$, where we get the picture
\begin{equation}
\xymatrix@R=5pt@C=40pt{%
  \Bbb R \ar@{<->}[dr] & \Bbb R^n & \La^2\Bbb R^n\\
  \Bbb R^n & \Bbb R\oplus S^2_0\Bbb R^n\oplus\La^2\Bbb R^n &
  \Bbb R^n\oplus
  W_2\oplus \La^3\Bbb R^n \ar@{<->}[ul] \\
  \Bbb R & \Bbb R^n \ar@{<->}[ul] & \La^2\Bbb R^n \ar@{<->}[ul] }
\label{example-decomp}\end{equation}
and we have indicated some components which are isomorphic as
representations of $O(n)$. Observe that assigning homogeneity $k+i$ to elements
of $\La^k\Bbb R^n\otimes\Bbb V_i$, we have chosen to identify
components of the same homogeneity.

We will make these identifications explicit in the language of bundles
immediately, but let us first state how we will use them. For the left
column, we will on the one hand define $\partial:\Bbb V\to \Bbb
R^n\otimes\Bbb V$, which vanishes on $\Bbb V_0$ and is injective on
$\Bbb V_1\oplus\Bbb V_2$. On the other hand, we will define
$\delta^*:\Bbb R^n\otimes\Bbb V\to\Bbb V$ by using inverse
identifications. For the right hand column, we will only use the
identifications from right to left to define $\delta^*:\La^2\Bbb
R^n\otimes\Bbb R\to\Bbb R^n\otimes\Bbb V$. Evidently, this map has
values in the kernel of $\delta^*:\Bbb R^n\otimes\Bbb V\to\Bbb V$, so
$\delta^*\o\delta^*=0$. By constructions, all these maps preserve
homogeneity. We also observe that $\ker(\delta^*)=S^2_0\Bbb
R^n\oplus\im(\delta^*)\subset\Bbb R^n\otimes\Bbb V$.

Now we can carry all this over to any Riemannian manifold of dimension
$n$. Sections of the bundle $V$ corresponding to $\Bbb V$ can be
viewed as triples consisting of two functions and a one--form. Since
the representation $\Bbb R^n$ corresponds to $T^*M$, the bundle
corresponding to $\La^k\Bbb R^n\otimes\Bbb V$ is $\La^kT^*M\otimes V$.
Sections of this bundle are triples consisting of two $k$--forms and
one $T^*M$--valued $k$--form. If there is no danger of confusion with
abstract indices, we will use subscripts $i=0,1,2$ to denote the
component of a section in $\La^kT^*M\otimes V_i$. To specify our maps,
we use abstract index notation and define $\partial:V\to T^*M\otimes
V$, $\delta^*:T^*M\otimes V\to V$ and $\delta^*:\La^2T^*M\otimes V\to
T^*M\otimes V$ by
$$
\partial\begin{pmatrix} h \\ \ph_b \\ f\end{pmatrix}:=\begin{pmatrix}
  0 \\ hg_{ab} \\ -\ph_a \end{pmatrix}\quad 
  \delta^*\begin{pmatrix} h_b \\ \ph_{bc} \\
        f_b\end{pmatrix}:= \begin{pmatrix} \tfrac{1}{n}\ph^c_c \\ -f_b
  \\ 0\end{pmatrix} \quad  \delta^*\begin{pmatrix} h_{ab} \\ \ph_{abc} \\
        f_{ab}\end{pmatrix}:= \begin{pmatrix}
  \tfrac{-1}{n-1}\ph_{ac}{}^c \\ \tfrac{1}{2}f_{ab}
  \\ 0\end{pmatrix}
$$
The numerical factors are chosen in such a way that our example
fits into the general framework developed in section \ref{3}.

We can differentiate sections of $V$ using the component--wise
Levi--Civita connection, which we denote by $\nabla$. Note that this
raises homogeneity by one. The core of the method is to mix this
differential term with an algebraic one. We consider the operation
$\Ga(V)\to\Om^1(M,V)$ defined by $\Si\mapsto \nabla\Si+\partial\Si$.
Since $\partial$ is tensorial and linear, this defines a linear
connection $\tilde\nabla$ on the vector bundle $V$.

We are ready to define the class of systems that we will look at.
Choose a bundle map (not necessarily linear) $A:V_0\oplus V_1\to
S^2_0T^*M$, and view it as $A:V\to T^*M\otimes V$. Notice that this
implies that $A$ increases homogeneities. Then consider the system
\begin{equation}
  \label{example-basic}
  \tilde\nabla\Si+A(\Si)=\delta^*\ps\qquad\text{for some\
  }\ps\in\Om^2(M,V). 
\end{equation}
We will show that on the one hand, this is equivalent to a second
order system on the $V_0$--component $\Si_0$ of $\Si$ and on the other
hand, it is equivalent to a first order system on $\Si$ in closed form.

\subsection{The first splitting operator}\label{2.2}
Since $A$ by definition has values in $\ker(\delta^*)$ and
$\delta^*\o\delta^*=0$, the system \eqref{example-basic} implies
$\delta^*(\tilde\nabla\Si)=0$. Hence we first have to analyze the
operator $\delta^*\o\tilde\nabla:\Ga(V)\to\Ga(V)$. Using abstract
indices and denoting the Levi--Civita connection by $\nabla_a$ we
obtain
$$
\Si=\begin{pmatrix} h \\ \ph_b \\ 
  f\end{pmatrix}\overset{\tilde\nabla_a}{\mapsto}
\begin{pmatrix}  \nabla_a h \\ \nabla_a\ph_b+hg_{ab} \\ 
  \nabla_a f-\ph_a\end{pmatrix}\overset{\delta^*}{\mapsto}
\begin{pmatrix}  \tfrac{1}{n}\nabla^b\ph_b+h \\ -\nabla_a f+\ph_a \\
  0\end{pmatrix}
$$
From this we can read off the set of all solutions of
$\delta^*\tilde\nabla\Si=0$. We can arbitrarily choose $f$. Vanishing
of the middle row then forces $\ph_a=\nabla_a f$, and inserting this,
vanishing of the top row is equivalent to
$h=-\tfrac{1}{n}\nabla^b\nabla_b f=-\tfrac{1}{n}\Delta f$, where
$\Delta$ denotes the Laplacian. Hence we get
\begin{proposition}
  For any $f\in C^\infty(M,\Bbb R)$, there is a unique $\Si\in\Ga(V)$
  such that $\Si_0=f$ and $\delta^*(\nabla\Si+\delta\Si)=0$. Mapping
  $f$ to this unique $\Si$ defines a second order linear differential
  operator $L:\Ga(V_0)\to\Ga(V)$, which is explicitly given by
  $$
  L(f)=\begin{pmatrix}-\tfrac{1}{n}\Delta f \\ \nabla_a f\\
    f\end{pmatrix}=\sum_{i=0}^2(-1)^i(\delta^*\nabla)^i
\begin{pmatrix} 0 \\ 0\\ f\end{pmatrix}.
$$ 
\end{proposition}

The natural interpretation of this result is that $V_0$ is viewed as a
quotient bundle of $V$, so we have the tensorial projection
$\Si\mapsto\Si_0$. The operator $L$ constructed provides a
differential splitting of this tensorial projection, which is
characterized by the simple property that its values are in the kernel
of $\delta^*\tilde\nabla$. Therefore, $L$ and its generalizations are
called \textit{splitting operators}.

\subsection{Rewriting as a higher order system}\label{2.3}
We have just seen that the system $\tilde\nabla\Si+A(\Si)=\delta^*\ps$
from \ref{2.1} implies that $\Si=L(f)$, where $f=\Si_0$. Now by
Proposition \ref{2.2}, the components of $L(f)$ in $V_0$ and $V_1$ are
$f$ and $\nabla f$, respectively. Hence $f\mapsto A(L(f))$ is a first
order differential operator $\Ga(V_0)\to\Ga(S^2_0T^*M)$. Conversely,
any first order operator
$D_1:\Ga(V_0)\to\Ga(S^2_0T^*M)\subset\Om^1(M,V)$ can be written as
$D_1(f)=A(L(f))$ for some $A:V\to T^*M\otimes V$ as in \ref{2.1}.

Next, for $f\in\Ga(V_0)$ we compute
$$
\tilde\nabla L(f)=\tilde\nabla_a \begin{pmatrix} -\tfrac{1}{n}\Delta f\\ \nabla_b f\\
  f\end{pmatrix}=
\begin{pmatrix}  -\tfrac{1}{n}\nabla_a\Delta f \\ \nabla_a\nabla_b
  f-\tfrac{1}{n}g_{ab}\Delta f\\ 0 \end{pmatrix}.
$$
The middle component of this expression is the tracefree part
$\nabla_{(a}\nabla_{b)_0}f$ of $\nabla^2 f$.

\begin{proposition}
  For any operator $D_1:C^\infty(M,\Bbb R)\to\Ga(S^2_0T^*M)$ of first
  order, there is a bundle map $A:V\to T^*M\otimes V$ such that
  $f\mapsto L(f)$ and $\Si\mapsto\Si_0$ induce inverse bijections
  between the sets of solutions of
\begin{equation}\label{example-ho}
\nabla_{(a}\nabla_{b)_0}f+D_1(f)=0
\end{equation} 
and of the basic system \eqref{example-basic}.
\end{proposition}
\begin{proof}
  We can choose $A:V_0\oplus V_1\to S^2_0T^*M\subset T^*M\otimes V$ in
  such a way that $D_1(f)=A(L(f))$ for all $f\in\Ga(V_0)$. From above
  we see that $\tilde\nabla L(f)+A(L(f))$ has vanishing bottom
  component and middle component equal to
  $\nabla_{(a}\nabla_{b)_0}f+D_1(f)$. From \eqref{example-decomp} we
  see that sections of $\im(\delta^*)\subset T^*M\otimes V$ are
  characterized by the facts that the bottom component vanishes, while
  the middle one is skew symmetric. Hence $L(f)$ solves
  \eqref{example-basic} if and only if $f$ solves \eqref{example-ho}.
  Conversely, we know from \ref{2.2} any solution $\Si$ of
  \eqref{example-basic} satisfies $\Si=L(\Si_0)$, and the result
  follows.
\end{proof}

Notice that in this result we do not require $D_1$ to be linear. In
technical terms, an operator can be written in the form
$f\mapsto\nabla_{(a}\nabla_{b)_0}f+D_1(f)$ for a first order operator
$D_1$, if and only if it is of second order, quasi--linear and its
principal symbol is the projection $S^2T^*M\to S^2_0T^*M$ onto the
tracefree part.

\subsection{Rewriting in closed form}\label{2.4}
Suppose that $\Si$ is a solution of \eqref{example-basic},
i.e.~$\tilde\nabla\Si+A(\Si)=\delta^*\ps$ for some $\ps$. Then the
discussion in \ref{2.3} shows that the two bottom components of
$\tilde\nabla\Si+A(\Si)$ actually have to vanish. Denoting the
components of $\Si$ as before, there must be a one--form $\tau_a$ such
that
\begin{equation}\label{almost-closed}
\begin{pmatrix}
  \nabla_ah+\tau_a\\ \nabla_a\ph_b+hg_{ab}+A_{ab}(f,\ph) \\
  \nabla_af-\ph_a
\end{pmatrix}=0. 
\end{equation}
Apart from the occurrence of $\tau_a$, this is a first order system in
closed form, so it remains to compute this one--form.

To do this, we use the \textit{covariant exterior derivative}
$d^{\tilde\nabla}:\Om^1(M,V)\to\Om^2(M,V)$ associated to
$\tilde\nabla$. This is obtained by coupling the exterior derivative
to the connection $\tilde\nabla$, so in particular on one--forms we
obtain
$$
d^{\tilde\nabla}\om(\xi,\eta)=\tilde\nabla_\xi(\om(\eta))-
\tilde\nabla_\eta(\om(\xi))-\om([\xi,\eta]). 
$$
Explicitly, on $\Om^1(M,V)$ the operator $d^{\tilde\nabla}$ is
given by
\begin{equation}\label{example-cov-ext}
\begin{pmatrix} h_b \\ \ph_{bc} \\
  f_b\end{pmatrix}\mapsto 2\begin{pmatrix} \nabla_{[a} h_{b]} \\ 
  \nabla_{[a}\ph_{b]c}-h_{[a}g_{b]c}\\ 
  \nabla_{[a}f_{b]}+\ph_{[ab]}\end{pmatrix},
\end{equation}
where square brackets indicate indicate an alternation of abstract
indices. 

Now almost by definition, $d^{\tilde\nabla}\tilde\nabla\Si$ is given
by the action of the curvature of $\tilde\nabla$ on $\Si$. One easily
computes directly that this coincides with the component--wise action
of the Riemann curvature. In particular, this is only non--trivial on
the middle component. On the other hand, since $A(\Si)=A_{ab}(f,\ph)$
is symmetric, we see that $d^{\tilde\nabla}(A(\Si))$ is concentrated
in the middle component, and it certainly can be written as
$\Ph_{abc}(f,\nabla f,\ph,\nabla\ph)$ for an appropriate bundle map
$\Ph$.  Together with the explicit formula, this shows that applying
the covariant exterior derivative to \eqref{almost-closed} we obtain
$$
\begin{pmatrix}
  2\nabla_{[a}\tau_{b]}\\ -R_{ab}{}^d{}_c\ph_d+\Ph_{abc}(f,\nabla
  f,\ph,\nabla\ph)-2\tau_{[a}g_{b]c}\\ 0
\end{pmatrix}=0.
$$
Applying $\delta^*$, we obtain an element with the bottom two rows
equal to zero and top row given by
$$
\tfrac{1}{n-1}\big(R_a{}^c{}^d{}_c\ph_d-\Ph_a{}^c{}_c(f,\nabla
f,\ph,\nabla\ph)\big)+\tau_a,
$$
which gives a formula for $\tau_a$. Finally, we define a bundle map
$C:V\to T^*M\otimes V$ by
$$
C\begin{pmatrix} h\\ \ph_b\\ f\end{pmatrix}:=
\begin{pmatrix}
\tfrac{-1}{n-1}\big(R_a{}^c{}^d{}_c\ph_d-
\Ph_a{}^c{}_c(f,\ph,\ph,-hg-A(f,\ph))\big)\\
  A_{ab}(f,\ph)\\ 0
\end{pmatrix}
$$ 
to obtain 
\begin{theorem}
  Let $D:C^\infty(M,\Bbb R)\to\Ga(S^2_0T^*M)$ be a quasi--linear
  differential operator of second order whose principal symbol is the
  projection $S^2T^*M\to S^2_0T^*M$ onto the tracefree part. Then
  there is a bundle map $C:V\to T^*M\otimes V$ which has the property
  that $f\mapsto L(f)$ and $\Si\mapsto\Si_0$ induce inverse bijections
  between the sets of solutions of $D(f)=0$ and of
  $\tilde\nabla\Si+C(\Si)=0$.  If $D$ is linear, then $C$ can be
  chosen to be a vector bundle map.
\end{theorem}

Since for any bundle map $C$, a solution of $\tilde\nabla\Si+C(\Si)=0$
is determined by its value in a single point, we conclude that any
solution of $D(f)=0$ is uniquely determined by the values of $f$,
$\nabla f$ and $\Delta f$ in one point.  Moreover, if $D$ is linear,
then the dimension of the space of solutions is always $\leq n+2$. In
this case, $\tilde\nabla+C$ defines a linear connection on the bundle
$V$, and the maximal dimension can be only attained if this connection
is flat.

Let us make the last step explicit for
$D(f)=\nabla_{(a}\nabla_{b)_0}f+A_{ab}f$ with some fixed section
$A_{ab}\in\Ga(S^2_0T^*M)$. From formula \eqref{example-cov-ext} we
conclude that
$$
\Ph_{abc}(f,\nabla
f,\ph,\nabla\ph)=2f\nabla_{[a}A_{b]c}+2A_{c[b}\nabla_{a]}f,
$$
and inserting we obtain the closed system
$$
\begin{cases}
  \nabla_ah-\tfrac{1}{n-1}\big(R_a{}^c{}^d{}_c\ph_d+
f\nabla^cA_{ac}+\ph^cA_{ac}\big)=0\\
\nabla_a\ph_b+hg_{ab}+fA_{ab}=0\\
\nabla_af-\ph_a=0 
\end{cases}
$$
which is equivalent to $\nabla_{(a}\nabla_{b)_0}f+A_{ab}f=0$.

\subsection{Remark}\label{2.5}
As a slight detour (which however is very useful for the purpose of
motivation) let me explain why the equation
$\nabla_{(a}\nabla_{b)_0}f+A_{ab}f=0$ is of geometric interest.  Let
us suppose that $f$ is a nonzero function.  The we can use it to
conformally rescale the metric $g$ to $\hat g:=\tfrac{1}{f^2} g$. Now
one can compute how a conformal rescaling affects various quantities,
for example the Levi--Civita connection. In particular, we can look at
the conformal behavior of the Riemannian curvature tensor. Recall that
the Riemann curvature can be decomposed into various components
according to the decomposition of $S^2(\La^2\Bbb R^n)$ as a
representation of $O(n)$. The highest weight part is the \textit{Weyl
  curvature}, which is independent of conformal rescalings.

Contracting the Riemann curvature via
$\text{Ric}_{ab}:=R_{ca}{}^c{}_b$, one obtains the \textit{Ricci
  curvature}, which is a symmetric two tensor. This can be further
decomposed into the \textit{scalar curvature} $R:=\text{Ric}^a{}_a$
and the tracefree part
$\text{Ric}^0_{ab}=\text{Ric}_{ab}-\tfrac{1}{n}Rg_{ab}$. Recall that a
Riemannian metric is called an \textit{Einstein metric} if the Ricci
curvature is proportional to the metric, i.e.~if
$\text{Ric}^0_{ab}=0$.

The behavior of the tracefree part of the Ricci curvature under a
conformal change $\hat g:=\tfrac{1}{f^2} g$ is easily determined
explicitly, see~\cite{BEG}. In particular, $\hat g$ is Einstein if and
only if
$$
\nabla_{(a}\nabla_{b)_0}f+A_{ij} f=0
$$
for an appropriately chosen $A_{ij}\in\Ga(S^2_0T^*M)$. Hence
existence of a nowhere vanishing solution to this equation is
equivalent to the possibility of rescaling $g$ conformally to an
Einstein metric. 

From above we know that for a general non--trivial solution $f$ of
this system and each $x\in M$, at least one of $f(x)$, $\nabla f(x)$,
and $\Delta f(x)$ must be nonzero. Hence $\{x:f(x)\neq 0\}$ is a dense
open subset of $M$, and one obtains a conformal rescaling to an
Einstein metric on this subset.

\section{The general procedure}\label{3}
\setcounter{lemma}1
\setcounter{theorem}1
\setcounter{proposition}2
The procedure carried out in an example in section \ref{2} can be
vastly generalized by replacing the standard representation by an
arbitrary irreducible representation of $G\cong O(n+1,1)$. (Things
also work for spinor representations, if one uses $Spin(n+1,1)$
instead.)  However, one has to replace direct observations by tools
from representation theory, and we discuss in this section, how this
is done.

\subsection{The Lie algebra $\frak{o}(n+1,1)$}\label{3.1}
We first have to look at the Lie algebra $\frak g\cong\frak{o}(n+1,1)$
of $G=O(\Bbb V)$. For the choice of inner product used in \ref{2.1}
this has the form
$$
\frak g=\left\{\begin{pmatrix} a & Z & 0 \\ X & A & -Z^t \\ 0 & -X^t &
    -a\end{pmatrix}: A\in\frak o(n), a\in\Bbb R, X\in\Bbb R^n,
  Z\in\Bbb R^{n*},\right\}. 
$$
The central block formed by $A$ represents the subgroup $O(n)$. The
element $E:=\left(\begin{smallmatrix} 1 & 0 & 0\\ 0 & 0 & 0\\ 0 & 0 &
    -1\end{smallmatrix}\right)$ is called the \textit{grading
  element}. Forming the commutator with $E$ is a diagonalizable map
$\frak g\to\frak g$ with eigenvalues $-1$, $0$, and $1$, and we denote
by $\frak g_i$ the eigenspace for the eigenvalue $i$. Hence $\fg_{-1}$
corresponds to $X$, $\fg_1$ to $Z$ and $\fg_0$ to $A$ and $a$.
Moreover, the Jacobi identity immediately implies that
$[\fg_i,\fg_j]\subset\fg_{i+j}$ with the convention that
$\fg_{i+j}=\{0\}$ unless $i+j\in \{-1,0,1\}$. Such a decomposition is
called a \textit{$|1|$--grading} of $\fg$. In particular, restricting
the adjoint action to $\frak o(n)$, one obtains actions on $\fg_{-1}$
and $\fg_1$, which are the standard representation respectively its
dual (and hence isomorphic to the standard representation).

Since the grading element $E$ acts diagonalizably under the adjoint
representation, it also acts diagonalizably on any finite dimensional
irreducible representation $\Bbb W$ of $\frak g$. If $w\in\Bbb W$ is
an eigenvector for the eigenvalue $j$, and $Y\in\fg_i$, then
$E\cdot Y\cdot w=Y\cdot E\cdot w+[E,Y]\cdot w$
shows that $Y\cdot w$ is an eigenvector with eigenvalue $i+j$. From
irreducibility it follows easily that denoting by $j_0$ the lowest
eigenvalue, the set of eigenvalues is $\{j_0,j_0+1,\dots,j_0+N\}$ for
some $N\geq 1$. Correspondingly, we obtain a decomposition $\Bbb
W=\Bbb W_0\oplus\dots\oplus \Bbb W_N$ such that $\fg_i\cdot\Bbb
W_j\subset\Bbb W_{i+j}$. In particular, each of the subspaces $\Bbb
W_j$ is invariant under the action of $\fg_0$ and hence in particular
under the action of $\frak o(n)$. Notice that the decomposition $\Bbb
V=\Bbb V_0\oplus\Bbb V_1\oplus\Bbb V_2$ used in section \ref{2} is
obtained in this way.

One can find a Cartan subalgebra of (the complexification of) $\frak
g$ which is spanned by $E$ and a Cartan subalgebra of (the
complexification of) $\frak o(n)$. The theorem of the highest weight
then leads to a bijective correspondence between finite dimensional
irreducible representations $\Bbb W$ of $\frak g$ and pairs $(\Bbb
W_0,r)$, where $\Bbb W_0$ is a finite dimensional irreducible
representation of $\frak o(n)$ and $r\geq 1$ is an integer. Basically,
the highest weight of $\Bbb W_0$ is the restriction to the Cartan
subalgebra of $\frak o(n)$ of the highest weight of $\Bbb W$, while
$r$ is related to the value of the highest weight on $E$. As the
notation suggests, we can arrange things in such a way that $\Bbb W_0$
is the lowest eigenspace for the action of $E$ on $\Bbb W$. For
example, the standard representation $\Bbb V$ in this notation
corresponds to $(\Bbb R,2)$. The explicit version of this
correspondence is not too important here, it is described in terms of
highest weights in \cite{BCEG} and in terms of Young diagrams in
\cite{Mike:prolon}. It turns out that, given $\Bbb W_0$ and $r$, the number
$N$ which describes the length of the grading can be easily computed.

\subsection{Kostant's version of the Bott--Borel--Weil theorem}\label{3.2} 
Suppose that $\Bbb W$ is a finite dimensional irreducible
representation of $\frak g$, decomposed as $\Bbb
W_0\oplus\dots\oplus\Bbb W_N$ as above. Then we can view $\La^k\Bbb
R^n\otimes\Bbb W$ as $\La^k\fg_1\otimes\Bbb W$, which leads to two
natural families of $O(n)$--equivariant maps. First we define
$\partial^*:\La^k\fg_1\otimes\Bbb W\to\La^{k-1}\fg_1\otimes\Bbb W$ by
$$
\partial^*(Z_1\wedge\dots\wedge Z_k\otimes w):=\textstyle\sum_{i=1}^k (-1)^i
Z_1\wedge\dots \wedge\widehat{Z_i}\wedge\dots\wedge Z_k\otimes
Z_i\cdot w,
$$
where the hat denotes omission. Note that if $w\in\Bbb W_j$, then
$Z_i\cdot w\in W_{j+1}$, so this operation preserves homogeneity. On
the other hand, we have $Z_i\cdot Z_j\cdot w-Z_j\cdot Z_i\cdot
w=[Z_i,Z_j]\cdot w=0$, since $\fg_1$ is a commutative subalgebra. This
easily implies that $\partial^*\o\partial^*=0$.

Next, there is an evident duality between $\fg_{-1}$ and $\fg_1$,
which is compatible with Lie theoretic methods since it is induced by
the Killing form of $\fg$. Using this, we can identify
$\La^k\fg_1\otimes\Bbb W$ with the space of $k$--linear alternating
maps $\fg_{-1}^k\to\Bbb W$. This gives rise to a natural map
$\partial=\partial_k:\La^k\fg_1\otimes\Bbb
W\to\La^{k+1}\fg_1\otimes\Bbb W$ defined by
$$
\partial\al(X_0,\dots,X_k):=
\textstyle\sum_{i=0}^k(-1)^iX_i\cdot\al(X_0,\dots,\widehat{X_i},\dots,X_k). 
$$
In this picture, homogeneity boils down to the usual notion for
multilinear maps, i.e.~$\al:(\fg_{-1})^k\to\Bbb W$ is homogeneous of
degree $\ell$ if it has values in $\Bbb W_{\ell-k}$. From this it
follows immediately that $\partial$ preserves homogeneities, and
$\partial\o\partial=0$ since $\fg_{-1}$ is commutative.

As a first step towards the proof of his version of the
Bott--Borel--Weil--theorem (see \cite{Kostant}), B.~Kostant proved the
following result:
\begin{lemma}
  The maps $\partial$ and $\partial^*$ are adjoint with respect to an
  inner product of Lie theoretic origin. For each degree $k$, one
  obtains an algebraic Hodge decomposition
$$
\La^k\fg_1\otimes\Bbb W=\im(\partial)\oplus
(\ker(\partial)\cap\ker(\partial^*))\oplus\im(\partial^*),
$$ 
with the first two summands adding up to $\ker(\partial)$ and the last
two summands adding up to $\ker(\partial^*)$. 

In particular, the restrictions of the canonical projections to the
subspace $\Bbb H_k:=\ker(\partial)\cap\ker(\partial^*)$ induce
isomorphisms $\Bbb H_k\cong\ker(\partial)/\im(\partial)$ and $\Bbb
H_k\cong\ker(\partial^*)/\im(\partial^*)$.
\end{lemma}

Since $\partial$ and $\partial^*$ are $\fg_0$--equivariant, all spaces
in the lemma are naturally representations of $\fg_0$ and all
statements include the $\fg_0$--module structure. Looking at the Hodge
decomposition more closely, we see that for each $k$, the map
$\partial$ induces an isomorphism 
$$
\La^k\fg_1\otimes\Bbb
W\supset\im(\partial^*)\to\im(\partial)\subset
\La^{k+1}\fg_1\otimes\Bbb W,
$$
while $\partial^*$ induces an isomorphism in the opposite
direction. In general, these two map are not inverse to each other, so
we replace $\partial^*$ by the map $\delta^*$ which vanishes on
$\ker(\partial^*)$ and is inverse to $\partial$ on $\im(\partial)$. Of
course, $\delta^*\o\delta^*=0$ and it computes the same cohomology as
$\partial^*$.

Kostant's version of the BBW--theorem computes (in a more general
setting to be discussed below) the representations $\Bbb H_k$ in an
explicit and algorithmic way. We only need the cases $k=0$ and $k=1$
here, but to formulate the result for $k=1$ we need a bit of
background. Suppose that $\Bbb E$ and $\Bbb F$ are finite dimensional
representations of a semisimple Lie algebra. Then the tensor product
$\Bbb E\otimes\Bbb F$ contains a unique irreducible component whose
highest weight is the sum of the highest weights of $\Bbb E$ and $\Bbb
F$. This component is called the \textit{Cartan product} of $\Bbb E$
and $\Bbb F$ and denoted by $\Bbb E\circledcirc \Bbb F$. Moreover,
there is a nonzero equivariant map $\Bbb E\otimes\Bbb F\to \Bbb
E\circledcirc\Bbb F$, which is unique up to multiples. This
equivariant map is also referred to as the \textit{Cartan product}.

The part of Kostant's version of the BBW--theorem that we need (proved
in \cite{BCEG} in this form) reads as follows,
\begin{theorem}
  Let $\Bbb W=\Bbb W_0\oplus\dots\oplus\Bbb W_N$ be the irreducible
  representation of $\frak g$ corresponding to the pair $(\Bbb
  W_0,r)$. 

\noindent
(i) In degree zero, $\im(\partial^*)=\Bbb W_1\oplus\dots\oplus\Bbb
W_N$ and $\Bbb H_0=\ker(\partial)=\Bbb W_0$.

\noindent
(ii) The subspace $\Bbb H_1\subset\fg_1\otimes\Bbb W$ is isomorphic to
$S^r_0\fg_1\circledcirc\Bbb W_0$. It is contained in $\fg_1\otimes\Bbb
W_{r-1}$ and it is the only irreducible component of
$\La^*\fg_1\otimes \Bbb W$ of this isomorphism type.
\end{theorem}

\subsection{Some more algebra}\label{3.2a}
Using Theorem \ref{3.2} we can now deduce the key algebraic ingredient
for the procedure. For each $i\geq 1$ we have $\partial:\Bbb
W_i\to\fg_1\otimes\Bbb W_{i-1}$. Next, we consider
$(\id\otimes\partial)\o\partial:\Bbb W_i\to\otimes^2\fg_1\otimes\Bbb
W_{i-2}$, and so on, to obtain $\fg_0$--equivariant maps
$$
\ph_i:=(\id\otimes\dots\otimes\id\otimes\partial)\o\dots\o
(\id\otimes\partial)\otimes\partial:\Bbb W_i\to\otimes^i\fg_1\otimes\Bbb W_0
$$ 
for $i=1,\dots,N$, and we put $\ph_0=\id_{\Bbb W_0}$.

\begin{proposition}
Let $\Bbb W=\Bbb W_0\oplus\dots\oplus\Bbb W_N$ correspond to $(\Bbb
W_0,r)$ and let $\Bbb K\subset S^r\fg_1\otimes\Bbb W_0$ be the kernel
of the Cartan product. Then we have

\noindent
(i) For each $i$, the map $\ph_i:\Bbb W_i\to\otimes^i\fg_1\otimes\Bbb
W_0$ is injective and hence an isomorphism onto its image. This image
is given by
$$
\im(\ph_i)=
\begin{cases}
  S^i\fg_1\otimes\Bbb W_0\qquad i<r\\
  (S^i\fg_1\otimes\Bbb W_0)\cap (S^{i-r}\fg_1\otimes\Bbb K)\qquad
  i\geq r.
\end{cases}
$$
(ii) For each $i<r$, the restriction of the map $\delta^*\otimes
\ph_{i-1}^{-1}$ to $S^i\fg_1\otimes\Bbb W_0\subset\fg_1\otimes
S^{i-1}\fg_1\otimes\Bbb W_0$ coincides with $\ph_i^{-1}$.
\end{proposition}
\begin{proof} (sketch) 
  (i) Part (i) of Theorem \ref{3.2} shows that $\partial:\Bbb
  W_i\to\fg_1\otimes\Bbb W_{i-1}$ is injective for each $i\geq 1$, so
  injectivity of the $\ph_i$ follows. Moreover, for $i\neq r$, the
  image of this map coincides with the kernel of
  $\partial_1:\fg_1\otimes\Bbb W_{i-1}\to \La^2\fg_1\otimes \Bbb
  W_{i-2}$, while for $i=r$ this kernel in addition contains a
  complementary subspace isomorphic to $S^k\fg_1\circledcirc\Bbb W_0$.
  A moment of thought shows that $\partial_1$ can be written as
  $2\text{Alt}\o(\id\otimes\partial_0)$, where $\text{Alt}$ denotes
  the alternation. This immediately implies that the $\ph_i$ all have
  values in $S^i\fg_1\otimes\Bbb W_0$ as well as the claim about the
  image for $i<r$. 

It further implies that $\id\otimes\ph_{r-1}$ restricts to isomorphisms
$$
\xymatrix@R=8pt{%
\fg_1\otimes\Bbb W_{r-1} \ar[r] & \fg_1\otimes S^{r-1}\fg_1\otimes\Bbb
W_0\\
\ker(\partial)\ar@{^{(}->}[u]\ar[r] & S^r \fg_1\otimes\Bbb
W_0\ar@{^{(}->}[u]\\
\im(\partial)\ar@{^{(}->}[u]\ar[r] & \Bbb K\ar@{^{(}->}[u],
}$$
which proves the claim on the image for $i=r$. For $i>r$ the claim
then follows easily as above. 

\noindent
(ii) This follows immediately from the fact that
$\delta^*|_{\im(\partial)}$ inverts $\partial|_{\im(\delta^*)}$.
\end{proof}

\subsection{Step one of the prolongation procedure}\label{3.3}
The developments in \ref{3.1}--\ref{3.2a} carry over to an arbitrary
Riemannian manifold $(M,g)$ of dimension $n$. The representation $\Bbb
W$ corresponds to a vector bundle $W=\oplus_{i=0}^NW_i$. Likewise,
$\Bbb H_1$ corresponds to a direct summand $H_1\subset T^*M\otimes
W_{r-1}$ which is isomorphic to $S^r_0T^*M\circledcirc W_0$. The maps
$\partial$, $\partial^*$, and $\delta^*$ induce vector bundle maps on
the bundles $\La^kT^*M\otimes W$ of $W$--valued differential forms,
and for $i=0,\dots,N$, the map $\ph_i$ induces a vector bundle map
$W_i\to S^iT^*M\otimes W_0$. We will denote all these maps by the same
symbols as their algebraic counterparts.  Finally, the Cartan product
gives rise to a vector bundle map $S^rT^*M\otimes W_0\to H_1$, which
is unique up to multiples.

We have the component--wise Levi--Civita connection $\nabla$ on $W$.
We will denote a typical section of $W$ by $\Si$. The subscript $i$
will indicate the component in $\La^kT^*M\otimes W_i$. Now we define a
linear connection $\tilde\nabla$ on $W$ by
$\tilde\nabla\Si:=\nabla\Si+\partial(\Si)$,
i.e.~$(\tilde\nabla\Si)_i=\nabla \Si_i+\partial(\Si_{i+1})$.  Next, we
choose a bundle map $A:W_0\oplus\dots\oplus W_{r-1}\to H_1$, view it
as $A:W\to T^*M\otimes W$ and consider the system
\begin{equation}
  \label{basic}
  \tilde\nabla \Si+A(\Si)=\delta^*\ps\quad\text{for some\
  }\ps\in\Om^2(M,W). 
\end{equation}
Since $A$ has values in $\ker(\delta^*)$ and $\delta^*\o\delta^*=0$,
any solution $\Si$ of this system has the property that
$\delta^*\tilde\nabla\Si=0$.

To rewrite the system equivalently as a higher order system, we define
a linear differential operator $L:\Ga(W_0)\to\Ga(W)$ by 
$$
L(f):=\textstyle\sum_{i=0}^N(-1)^i(\delta^*\o\nabla)^i f.
$$

\begin{proposition}
  (i) For $f\in\Ga(W_0)$ we have $L(f)_0=f$ and
  $\partial^*\tilde\nabla L(f)=0$, and $L(f)$ is uniquely determined
  by these two properties.

\noindent
(ii) For $\ell=0,\dots,N$ the component $L(f)_\ell$ depends only on
the $\ell$--jet of $f$. More precisely, denoting by $J^\ell W_0$ the
$\ell$th jet prolongation of the bundle $W_0$, the operator $L$
induces vector bundle maps $J^\ell W_0\to W_0\oplus\dots\oplus
W_\ell$, which are isomorphisms for all $\ell<r$.
\end{proposition}
\begin{proof}
  (i) Putting $\Si=L(f)$ it is evident that $\Si_0=f$ and
  $\Si_{i+1}=-\delta^*\nabla \Si_i$ for all $i\geq 0$. Therefore,
$$
(\tilde\nabla\Si)_i=\nabla\Si_i+\partial(\Si_{i+1})=\nabla\Si_i-
\partial\delta^*\nabla \Si_i
$$
for all $i$. Since $\delta^*\partial$ is the identity on
$\im(\delta^*)$, we get $\delta^*\tilde\nabla\Si=0$. 

Conversely, expanding the equation $0=\delta^*\tilde\nabla \Si$ in
components we obtain
$$
\Si_{i+1}=\delta^*\partial\Si_{i+1}=-\delta^*\nabla\Si_i,
$$ 
which inductively implies $\Si=L(\Si_0)$. 

\noindent
(ii) By definition, $L(f)_\ell$ depends only on $\ell$ derivatives of
$f$. Again by definition, $L(f)_1=\delta^*\nabla f$, and if $r>1$,
this equals $\ph_1^{-1}(\nabla f)$. Naturality of $\delta^*$ implies
that 
$$
L(f)_2=\delta^*\nabla\delta^*\nabla
f=\delta^*\o(\id\otimes\delta^*)(\nabla^2
f)=\delta^*\o(\id\otimes\ph_1^{-1})(\nabla^2 f).
$$
Replacing $\nabla^2$ by by its symmetrization changes the
expression by a term of order zero, so we see that, if $r>2$ and up to
lower order terms, $L(f)_2$ is obtained by applying $\ph_2^{-1}$ to
the symmetrization of $\nabla^2 f$. Using part (ii) of Proposition
\ref{3.2a} and induction, we conclude that for $\ell<r$ and up to
lower order terms $L(f)_\ell$ is obtained by applying $\ph_\ell^{-1}$
to the symmetrized $\ell$th covariant derivative of $\ell$, and the
claim follows.
\end{proof}

Note that part (ii) immediately implies that for a bundle map $A$ as
defined above, $f\mapsto A(L(f))$ is a differential operator
$\Ga(W_0)\to\Ga(H_1)$ of order at most $r-1$ and any such operator is
obtained in this way.

\subsection{The second step of the procedure}\label{3.4}
For a section $f\in\Ga(W_0)$ we next define $D^{\Bbb W}(f)\in\Ga(H_1)$
to be the component of $\tilde\nabla L(f)$ in
$\Ga(H_1)\subset\Om^1(M,W)$. We know that $(\tilde\nabla
L(f))_{r-1}=\nabla L(f)_{r-1}+\partial L(f)_r$, and the second summand
does not contribute to the $H_1$--component. Moreover, from the proof
of Proposition \ref{3.3} we know that, up to lower order terms,
$L(f)_{r-1}$ is obtained by applying $\ph_{r-1}^{-1}$ to the
symmetrized $(r-1)$--fold covariant derivative of $f$. Hence up to
lower order terms, $\nabla L(f)_{r-1}$ is obtained by applying
$\id\otimes\ph_{r-1}^{-1}$ to the symmetrized $r$--fold covariant
derivative of $f$. Using the proof of Proposition \ref{3.2a} this
easily implies that the principal symbol of $ D^{\Bbb W}$ is (a
nonzero multiple of) the Cartan product $S^rT^*M\otimes W_0\to
S^r_0T^*M\circledcirc W_0=H_1$.

\begin{proposition}
  Let $D:\Ga(W_0)\to\Ga(H_1)$ be a quasi--linear differential operator
  of order $r$ whose principal symbol is given by the Cartan product
  $S^rT^*M\otimes W_0\to S^r_0T^*M\circledcirc W_0$. Then there is a
  bundle map $A:W\to T^*M\otimes W$ as in \ref{3.3} such that
  $\Si\mapsto\Si_0$ and $f\mapsto L(f)$ induce inverse bijections
  between the sets of solutions of $D(f)=0$ and of the basic system
  \eqref{basic}.
\end{proposition}
\begin{proof}
  This is completely parallel to the proof of Proposition \ref{2.3}:
  The conditions on $D$ exaclty means that it can be written in the
  form $D(f)=D^{\Bbb W}(f)+A(L(f))$ for an appropriate choice of $A$
  as above. Then $\tilde\nabla L(f)+A(L(f))$ is a section of the
  subbundle $\ker(\delta^*)$ and the component in $H_1$ of this
  section equals $D(f)$. Of course, being a section of $\im(\delta^*)$
  is equivalent to vanishing of the $H_1$--component.
  
  Conversely, Proposition \ref{3.3} shows that any solution $\Si$ of
  \eqref{basic} is of the form $\Si=L(\Si_0)$.
\end{proof}

\subsection{The last step of the procedure}\label{3.5}
\setcounter{theorem}5
To rewrite the basic system \eqref{basic} in first order closed form,
we use the covariant exterior derivative $d^{\tilde\nabla}$. Suppose
that $\al\in\Om^1(M,W)$ has the property that its components $\al_i$
vanish for $i=0,\dots,\ell$. Then one immediately verifies that
$(d^{\tilde\nabla}\al)_i=0$ for $i=0,\dots,\ell-1$ and
$(d^{\tilde\nabla}\al)_\ell=\partial(\al_{\ell+1})$, so
$(\delta^*d^{\tilde\nabla}\al)_i$ vanishes for $i\leq\ell$ and equals
$\delta^*\partial\al_{\ell+1}$ for $i=\ell+1$. If we in addition
assume that $\al$ is a section of the subbundle $\im(\delta^*)$, then
the same is true for $\al_{\ell+1}$ and hence
$\delta^*\partial\al_{\ell+1}= \al_{\ell+1}$.

Suppose that $\Si$ solves the basic system \eqref{basic}. Then
applying $\delta^*d^{\tilde\nabla}$, we obtain
$$
\delta^*(R\bullet\Si
+d^{\tilde\nabla}(A(\Si)))=\delta^*d^{\tilde\nabla}\delta^*\ps,
$$
where we have used that, as in \ref{2.4},
$d^{\tilde\nabla}\tilde\nabla\Si$ is given by the action of the
Riemann curvature $R$. From above we see that we can compute the
lowest nonzero homogeneous component of $\delta^*\ps$ from this
equation. We can then move this to the other side in \eqref{basic} to 
obtain an equivalent system whose right hand side starts one
homogeneity higher. The lowest nonzero homogeneous component of the
right hand side can then be computed in the same way, and iterating
this we conclude that \eqref{basic} can be equivalently written as
\begin{equation}
  \label{hoe}
\tilde\nabla \Si+B(\Si)=0  
\end{equation}
for a certain differential operator $B:\Ga(W)\to\Om^1(M,W)$.

While $B$ is a higher order differential operator in general, it is
crucial that the construction gives us a precise control on the order
of the individual components of $B$. From the construction it follows
that $B(\Si)_i\in\Om^1(M,V_i)$ depends only on $\Si_0,\dots,\Si_i$,
and the dependence is tensorial in $\Si_i$, first order in $\Si_{i-1}$
and so on up to $i$th order in $\Si_0$.

In particular, the component of \eqref{hoe} in $\Om^1(M,W_0)$ has the
form $\nabla\Si_0=C_0(\Si_0,\Si_1)$. Next, the component in
$\Om^1(M,W_1)$ has the form $\nabla\Si_1=\tilde
C_1(\Si_0,\Si_1,\Si_2,\nabla\Si_0)$, and we define 
$$
C_1(\Si_0,\Si_1,\Si_2):=\tilde C_1(\Si_0,\Si_1,\Si_2,-C_0(\Si_0,\Si_1)).
$$ 
Hence the two lowest components of \eqref{hoe} are equivalent to 
$$
\begin{cases}
 \nabla\Si_1= C_1(\Si_0,\Si_1,\Si_2)\\
 \nabla\Si_0=C_0(\Si_0,\Si_1)
\end{cases}
$$
Differentiating the lower row and inserting for $\nabla\Si_0$ and
$\nabla\Si_1$ we get an expression for $\nabla^2\Si_0$ in terms of
$\Si_0,\Si_1,\Si_2$. Continuing in this way, one proves
\begin{theorem}
  Let $D:\Ga(W_0)\to\Ga(H_1)$ be a quasi--linear differential operator
  of order $r$ with principal symbol the Cartan product
  $S^rT^*M\otimes W_0\to S^rT^*M_0\circledcirc W_0$. Then there is a
  bundle map $C:W\to T^*M\otimes W$ such that $\Si\mapsto\Si_0$ and
  $f\mapsto L(f)$ induce inverse bijections between the sets of
  solutions of $D(f)=0$ and of $\tilde\nabla \Si+C(\Si)=0$. If $D$ is
  linear, then $C$ can be chosen to be a vector bundle map.
\end{theorem}

This in particular shows that any solution of $D(f)=0$ is determined
by the value of $L(f)$ in one point, end hence by the $N$--jet of $f$
in one point. For linear $D$, the dimension of the space of solutions
is bounded by $\dim(\Bbb W)$ and equality can be only attained if the
linear connection $\tilde\nabla+C$ on $W$ is flat. A crucial point
here is of course that $\Bbb W$, and hence $\dim(\Bbb W)$ and $N$ can
be immediately computed from $\Bbb W_0$ and $r$, so all this
information is available in advance, without going through the
procedure. As we shall see later, both the bound on the order and the
bound on the dimension are sharp.

To get a feeling for what is going on, let us consider some examples.
If we look at operators on smooth functions, we have $\Bbb W_0=\Bbb
R$. The representation associated to $(\Bbb R,r)$ is $S^{r-1}_0\Bbb
V$, the tracefree part of the $(r-1)$st symmetric power of the
standard representation $\Bbb V$. A moment of thought shows that the
eigenvalues of the grading element $E$ on this representation range
from $-r+1$ to $r-1$, so $N=2(r-1)$. On the other hand, for $r\geq 3$
we have
$$
\dim(S^{r-1}_0\Bbb V)=\dim(S^{r-1}\Bbb V)-\dim(S^{r-3}\Bbb
V)=(n+2r-2)\frac{(n+2r-2)!}{n!(r-1)!},
$$
and this is the maximal dimension of the space of solutions of any
system with principal part
$f\mapsto\nabla_{(a_1}\nabla_{a_2}\dots\nabla_{a_r)_0}f$ for $f\in
C^\infty(M,\Bbb R)$.

As an extreme example let us consider the conformal Killing equation
on tracefree symmetric tensors. Here $W_0=S^k_0TM$ for some $k$ and
$r=1$. The principal part in this case is simply
$$
f^{a_1\dots a_k}\mapsto\nabla^{(a}f^{a_1\dots a_k)_0}. 
$$
The relevant representation $\Bbb W$ in this case turns out to be
$\circledcirc^k\frak g$, i.e.~the highest weight subspace in $S^k\frak
g$. In particular $N=2k$ in this case, so even though we consider
first order systems, many derivatives a needed to pin down a solution.
The expression for $\dim(\Bbb W)$ is already reasonably complicated in
this case, namely (see \cite{Mike:symmetries})
$$
\dim(\Bbb W)=\frac{(n+k-3)!(n+k-2)!(n+2k)!}{k!(k+1)!(n-2)!n!(n+2k-3)!}
$$
The conformal Killing equation $\nabla^{(a}f^{a_1\dots a_k)_0}=0$
plays an important role in the description of symmetries of the
Laplacian on a Riemannian manifold, see \cite{Mike:symmetries}.

\section{Conformally invariant differential operators}\label{4}
We now move to the method for constructing conformally invariant
differential operators, which gave rise to the prolongation procedure
discussed in the last two sections.

\setcounter{proposition}3

\subsection{Conformal geometry}\label{4.1}
Let $M$ be a smooth manifold of dimension $n\geq 3$. As already
indicated in \ref{2.4}, two Riemannian metrics $g$ and $\hat g$ on $M$
are called \textit{conformally equivalent} if and only if there is a
positive smooth function $\ph$ on $M$ such that $\hat g=\ph^2 g$. A
\textit{conformal structure} on $M$ is a conformal equivalence class
$[g]$ of metrics, and then $(M,[g])$ is called a conformal manifold. A
\textit{conformal isometry} between conformal manifolds $(M,[g])$ and
$(\tilde M,[\tilde g])$ is a local diffeomorphism which pulls back one
(or equivalently any) metric from the class $[\tilde g]$ to a metric
in $[g]$.

A Riemannian metric on $M$ can be viewed as a reduction of structure
group of the frame bundle to $O(n)\subset GL(n,\Bbb R)$. In the same
way, a conformal structure is a reduction of structure group to
$CO(n)\subset GL(n,\Bbb R)$, the subgroup generated by $O(n)$ and
multiples of the identity. 

We want to clarify how the inclusion $O(n)\hookrightarrow G\cong
O(n+1,1)$ which was the basis for our prolongation procedure is
related to conformal geometry. For the basis $\{e_0,\dots,e_{n+1}\}$
used in \ref{2.1}, this inclusion was simply given by $A\mapsto
\left(\begin{smallmatrix} 1& 0 & 0 \\ 0 & A & 0\\ 0 & 0 &
    1\end{smallmatrix} \right)$. In \ref{3.1} we met the decomposition
$\fg=\fg_{-1}\oplus\fg_0\oplus\fg_1$ of the Lie algebra $\fg$ of $G$.
We observed that this decomposition is preserved by $O(n)\subset G$
and in that way $\fg_{\pm 1}$ is identified with the standard
representation. But there is a larger subgroup with these properties.
Namely, for elements of 
$$
G_0:=\left\{\left(\begin{smallmatrix} a & 0 & 0 \\ 0 & A & 0\\ 0 &
      0 & a^{-1}\end{smallmatrix} \right):a\in\Bbb R\setminus 0, A\in
  O(n)\right\}\subset G,
$$
the adjoint action preserves the grading, and maps $X\in\fg_{-1}$
to $a^{-1}AX$, so $G_0\cong CO(\fg_{-1})$. Note that $G_0\subset G$
corresponds to the Lie subalgebra $\fg_0\subset\fg$. 

Now there is a more conceptual way to understand this. Consider the
subalgebra $\frak p:=\fg_0\oplus\fg_1\subset\fg$ and let $P\subset G$
be the corresponding Lie subgroup. Then $P$ is the subgroup of
matrices which are block--upper--triangular with blocks of sizes $1$,
$n$, and $1$. Equivalently, $P$ is the stabilizer in $G$ of the
isotropic line spanned by the basis vector $e_0$. The group $G$ acts
transitively on the space of all isotropic lines in $\Bbb V$, so one
may identify this space with the homogeneous space $G/P$.

Taking coordinates $z_i$ with respect to an orthonormal basis of $\Bbb
V$ for which the first $n+1$ vectors are positive and the last one is
negative, a vector is isotropic if and only if
$\sum_{i=0}^nz_i^2=z_{n+1}^2$. Hence for a nonzero isotropic vector
the last coordinate is nonzero and any isotropic line contains a
unique vector whose last coordinate equals $1$. But this shows that the
space of isotropic lines in $\Bbb V$ is an $n$--sphere, so $G/P\cong
S^n$. 

Given a point $x\in G/P$, choosing a point $v$ in the corresponding
line gives rise to an identification $T_xS^n\cong v^\perp/\Bbb R v$
and that space carries a positive definite inner product induced by
$\langle\ ,\ \rangle$. Passing from $v$ to $\la v$, this inner product
gets scaled by $\la^2$, so we get a canonical conformal class of inner
products on each tangent space, i.e.~a conformal structure on $S^n$.
This conformal structure contains the round metric of $S^n$.

The action $\ell_g$ of $g\in G$ on the space of null lines by
construction preserves this conformal structure, so $G$ acts by
conformal isometries. It turns out, that this identifies
$G/\{\pm\id\}$ with the group of all conformal isometries of $S^n$.
For the base point $o=eP\in G/P$, the tangent space $T_o(G/P)$ is
naturally identified with $\frak g/\frak p\cong\fg_{-1}$. Let
$P_+\subset P$ be the subgroup of those $g\in P$ for which
$T_o\ell_g=\id$. Then one easily shows that $P/P_+\cong G_0$ and the
isomorphism $G_0\cong CO(\fg_{-1})$ is induced by $g\mapsto
T_o\ell_g$. Moreover, $P_+$ has Lie algebra $\frak g_1$ and
$\exp:\fg_1\to P_+$ is a diffeomorphism.

\subsection{Conformally invariant differential operators}\label{4.2} 
Let $(M,[g])$ be a conformal manifold. Choosing a metric $g$ from the
conformal class, we get the Levi--Civita connection $\nabla$ on each
Riemannian natural bundle as well as the Riemann curvature tensor $R$.
Using $g$, its inverse, and $R$, we can write down differential
operators, and see how they change if $g$ is replaced by a conformally
equivalent metric $\hat g$. Operators obtained in that way, which do
not change at all under conformal rescalings are called
\textit{conformally invariant}. In order to do this successfully one
either has to allow density bundles or deal with conformal weights,
but I will not go into these details here. The best known example of
such an operator is the conformal Laplacian or Yamabe operator which
is obtained by adding an appropriate amount of scalar curvature to the
standard Laplacian.

The definition of conformally invariant operators immediately suggests
a naive approach to their construction. First choose a principal part
for the operator. Then see how this behaves under conformal rescalings
and try to compensate the changes by adding lower order terms
involving curvature quantities. This approach (together with a bit of
representation theory) easily leads to a complete classification of
conformally invariant first order operators, see \cite{Fegan}. Passing
to higher orders, the direct methods get surprisingly quickly out of
hand.

The basis for more invariant approaches is provided by a classical
result of Elie Cartan, which interprets general conformal structures
as analogs of the homogeneous space $S^n\cong G/P$ from \ref{4.1}. As
we have noted above, a conformal structure $[g]$ on $M$ can be
interpreted as a reduction of structure group. This means that a
conformal manifold $(M,[g])$ naturally carries a principal bundle with
structure group $CO(n)$, the \textit{conformal frame bundle}. Recall
from \ref{4.1} that the conformal group $G_0=CO(\fg_{-1})\cong CO(n)$
can be naturally viewed as a quotient of the group $P$. Cartan's
result says that the conformal frame bundle can be canonically
extended to a principal fiber bundle $\Cal G\to M$ with structure
group $P$, and $\Cal G$ can be endowed with a canonical Cartan
connection $\om\in\Om^1(\Cal G,\frak g)$. The form $\om$ has similar
formal properties as the Maurer--Cartan form on $G$, i.e.~it defines a
trivialization of the tangent bundle $T\Cal G$, which is
$P$--equivariant and reproduces the generators of fundamental vector
fields.

While the canonical Cartan connection is conformally invariant, it is
not immediately clear how to use it to construct differential
operators. The problem is that, unlike principal connections, Cartan
connections do not induce linear connections on associated vector
bundles. 

\subsection{The setup for the conformal BGG machinery}\label{4.3}
Let us see how the basic developments from \ref{3.1}--\ref{3.2a}
comply with our new point of view. First of all, for $g\in P$, the
adjoint action does not preserve the grading
$\fg=\fg_{-1}\oplus\fg_0\oplus \fg_1$, but it preserves the
subalgebras $\frak p=\fg_0\oplus\fg_1$, and $\fg_1$. More generally,
if $\Bbb W=\Bbb W_0\oplus\dots\oplus\Bbb W_N$ is an irreducible
representation of $\fg$ decomposed according to eigenspaces of the
grading element $E$, then each of the subspaces $\Bbb
W_i\oplus\dots\oplus\Bbb W_N$ is $P$--invariant. Since $P$ naturally
acts on $\fg_1$ and on $\Bbb W$, we get induced actions on
$\La^k\fg_1\otimes\Bbb W$ for all $k$. The formula for
$\partial^*:\La^k\fg_1\otimes\Bbb W\to \La^{k-1}\fg_1\otimes\Bbb W$
uses only the action of $\fg_1$ on $\Bbb W$, so $\partial^*$ is
$P$--equivariant.

In contrast to this, the only way to make $P$ act on $\fg_{-1}$ is via
the identification with $\frak g/\frak p$. However, the action of
$\fg_{-1}$ on $\Bbb W$ has no natural interpretation in this
identification, and $\partial:\La^k\fg_1\otimes\Bbb W\to
\La^{k+1}\fg_1\otimes\Bbb W$ is \textit{not} $P$--equivariant.

Anyway, given a conformal manifold $(M,[g])$ we can now do the
following. Rather than viewing $\Bbb W$ just as sum of representations
of $G_0\cong CO(n)$, we can view it as a representation of $P$, and
form the associated bundle $\Cal W:=\Cal G\x_P\Bbb W\to M$. Bundles
obtained in this way are called \textit{tractor bundles}. I want to
emphasize at this point that the bundle $\Cal W$ is of completely
different nature than the bundle $W$ used in section \ref{3}. To see
this, recall that elements of the subgroup $P_+\subset P$ act on $G/P$
by diffeomorphisms which fix the base point $o=eP$ to first
order. Therefore, the action of such a diffeomorphism on the fiber
over $o$ of any tensor bundle is the identity. On the other hand, it
is easy to see that on the fiber over $o$ of any tractor bundle, this
action is always non--trivial. Hence tractor bundles are unusual
geometric objects. 

Examples of tractor bundles have already been introduced as an
alternative to Cartan's approach in the 1920's and 30's, in particular
in the work of Tracy Thomas, see \cite{Thomas}. Their key feature is
that the canonical Cartan connection $\om$ induces a canonical linear
connection, called the \textit{normal tractor connection} on each
tractor bundle. This is due to the fact that these bundles do not
correspond to general representations of $P$, but only to
representations which extend to the big group $G$. We will denote the
normal tractor connection on $\Cal W$ by $\nabla^{\Cal W}$. These
connections automatically combine algebraic and differential parts.

The duality between $\fg_1$ and $\fg_{-1}$ induced by the Killing
form, is more naturally viewed as a duality between $\fg_1$ and $\frak
g/\frak p$. Via the Cartan connection $\om$, the associated bundle
$\Cal G\x_P(\frak g/\frak p)$ is isomorphic to the tangent bundle
$TM$.  Thus, the bundle $\Cal G\x_P(\La^k\fg_1\otimes\Bbb W)$ is again
the bundle $\La^kT^*M\otimes\Cal W$ of $\Cal W$--valued forms. Now it
turns out that in a well defined sense (which however is rather awkward
to express), the lowest nonzero homogeneous component of $\nabla^{\Cal
  W}$ is of degree zero, it is tensorial and induced by the Lie
algebra differential $\partial$.

Equivariancy of $\partial^*$ implies that it defines bundle maps
$$
\partial^*:\La^kT^*M\otimes\Cal W\to\La^{k-1}T^*M\otimes\Cal W
$$
for each $k$. In particular,
$\im(\partial^*)\subset\ker(\partial^*)\subset \La^kT^*M\otimes\Cal W$
are natural subbundles, and we can form the subquotient
$H_k:=\ker(\partial^*)/\im(\partial^*)$. It turns out that these
bundles are always naturally associated to the conformal frame bundle,
so they are usual geometric objects like tensor bundles. The explicit
form of the bundles $H_k$ can be computed algorithmically using
Kostant's version of the BBW--theorem.

\subsection{The conformal BGG machinery}\label{4.4}
The normal tractor connection $\nabla^{\Cal W}$ extends to the
covariant exterior derivative, which we denote by $d^{\Cal
  W}:\Om^k(M,\Cal W)\to\Om^{k+1}(M,\Cal W)$. The lowest nonzero
homogeneous component of $d^{\Cal W}$ is of degree zero, tensorial,
and induced by $\partial$.

Now for each $k$, the operator $\partial^* d^{\Cal W}$ on $\Om^k(M,V)$
is conformally invariant and its lowest nonzero homogeneous component
is the tensorial map induced by $\partial^*\partial$. By Theorem
\ref{3.2}, $\partial^*\partial$ acts invertibly on
$\im(\partial^*)$. Hence we can find a (non--natural) bundle map $\be$
on $\im(\partial^*)$ such that $\be\partial^*d^\Cal W$ reproduces the
lowest nonzero homogeneous component of sections of
$\im(\partial^*)$. Therefore, the operator $\id-\be\partial^*d^\Cal
W$ is (at most $N$--step) nilpotent on $\Ga(\im(\partial^*))$, which
easily implies that 
$$
\left(\textstyle\sum_{i=0}^N(\id-\be\partial^*d^{\Cal W})^i\right)\be
$$
defines a differential operator $Q$ on $\Ga(\im(\partial^*))$ which
is inverse to $\partial^*d^{\Cal W}$ and therefore conformally
invariant.

Next, we have a canonical bundle map
$$
\pi_H:\ker(\partial^*)\to\ker(\partial^*)/\im(\partial^*)=H_k,
$$
and we denote by the same symbol the induced tensorial projection
on sections.  Given $f\in\Ga(H_k)$ we can choose $\ph\in\Om^k(M,\Cal
W)$ such that $\partial^*\ph=0$ and $\pi_H(\ph)=f$, and consider
$\ph-Q\partial^*d^{\Cal W}\ph$. By construction, $\ph$ is uniquely
determined up to adding sections of $\im(\partial^*)$. Since these are
reproduced by $Q\partial^*d^{\Cal W}$, the above element is
independent of the choice of $\ph$ and hence defines
$L(f)\in\Om^k(M,\Cal W)$. Since $Q$ has values in
$\Ga(\im(\partial^*))$ we see that $\pi_H(L(f))=f$, and since
$\partial^*d^{\Cal W}Q$ is the identity on $\Ga(\im(\partial^*))$ we
get $\partial^*d^{\Cal W}L(f)=0$. If $\ph$ satisfies $\pi_H(\ph)=f$
and $\partial^*d^{\Cal W}\ph=0$, then
$$
L(f)=\ph-Q\partial^*d^{\Cal W}\ph=\ph,
$$
so $L(f)$ is uniquely determined by these two properties.

By construction, the operator $L:\Ga(H_k)\to\Om^k(M,\Cal W)$ in
conformally invariant. Moreover, $d^{\Cal W}L(f)$ is a section of
$\ker(\partial^*)$, so we can finally define the BGG--operators
$D^{\Cal W}:\Ga(H_k)\to\Ga(H_{k+1})$ by $D^{\Cal W}(f):=\pi_H(d^{\Cal
  W}L(f))$. They are conformally invariant by construction.

To obtain additional information, we have to look at structures which
are locally conformally flat or equivalently locally conformally
isometric to the sphere $S^n$. It is a classical result that local
conformal flatness is equivalent to vanishing of the curvature of the
canonical Cartan connection. 

\begin{proposition}
  On locally conformally flat manifolds, the BGG operators form a
  complex $(\Ga(H_*),D^{\Cal W})$, which is a fine resolution of the
  constant sheaf $\Bbb W$. 
\end{proposition}
\begin{proof}
  The curvature of any tractor connection is induced by the Cartan
  curvature (see \cite{TAMS}), so on locally conformally flat
  structures, all tractor connections are flat. This implies that the
  covariant exterior derivative satisfies $d^{\Cal W}\o d^{\Cal W}=0$.
  Thus $(\Om^*(M,\Cal W),d^{\Cal W})$ is a fine resolution of the
  constant sheaf $\Bbb W$. 
  
  For $f\in\Ga(H_k)$ consider $d^{\Cal W}L(f)$. By construction, this
  lies in the kernel of $\partial^*$ and since $d^{\Cal W}\o d^{\Cal
    W}=0$, it also lies in the kernel of $\partial^*d^{\Cal W}$. From
  above we know that this implies that
$$
d^{\Cal W}L(f)=L(\pi_H(d^{\Cal W}L(f)))=L(D^{\Cal W}(f)).
$$
This shows that $L\o D^{\Cal W}\o D^{\Cal W}=d^{\Cal W}\o d^{\Cal
  W}\o L=0$ and hence $D^{\Cal W}\o D^{\Cal W}=0$, so
$(\Ga(H_*),D^{\Cal W})$ is a complex.  The operators $L$ define a
chain map from this complex to $(\Om^*(M,\Cal W),d^{\Cal W})$, and we
claim that this chain map induces an isomorphism in cohomology. 

First, take $\ph\in\Om^k(M,\Cal W)$ such that $d^{\Cal W}\ph=0$. Then
$$
\tilde \ph:=\ph-d^{\Cal W}Q\partial^*\ph
$$ 
is cohomologous to $\ph$ and satisfies $\partial^*\tilde\ph=0$.
Moreover, $d^{\Cal W}\tilde \ph=d^{\Cal W}\ph=0$, so $\partial^*
d^{\Cal W}\ph=0$. Hence $\tilde \ph=L(\pi_H(\tilde\ph))$ and $D^{\Cal
  W}(\pi_H(\tilde\ph))=0$, so the induced map in cohomology is
surjective.
  
Conversely, assume that $f\in\Ga(H_k)$ satisfies $D^{\Cal W}(f)=0$ and
that $L(f)=d^{\Cal W}\ph$ for some $\ph\in\Om^{k-1}(M,\Cal W)$.  As
before, replacing $\ph$ by $\ph-d^{\Cal W}Q\partial^*\ph$ we may
assume that $\partial^*\ph=0$. But together with $\partial^*L(f)=0$
this implies $\ph=L(\pi_H(\ph))$ and thus $f=\pi_H(L(f))=D^{\Cal
  W}(\pi_H(\ph))$. Hence the induced map in cohomology is injective,
too. Since this holds both locally and globally, the proof is
complete. 
\end{proof}

Via a duality between invariant differential operators and
homomorphisms of generalized Verma modules, this reproduces the
original BGG resolutions as constructed in \cite{Lepowsky}. Via the
classification of such homomorphisms, one also concludes that this
construction produces a large subclass of all those conformally
invariant operators which are non--trivial on locally conformally flat
structures.

Local exactness of the BGG sequence implies that all the operators
$D^{\Cal W}$ are nonzero on locally conformally flat manifolds.
Passing to general conformal structures does not change the principal
symbol of the operator $D^{\Cal W}$, so we always get non--trivial
operators.

On the other hand, we can also conclude that the bounds obtained from
Theorem \ref{3.5} are sharp. From Theorem \ref{3.2} we conclude that
any choice of metric in the conformal class identifies $H_0=\Cal
W/\im(\partial^*)$ with the bundle $W_0$ and $H_1$ with its
counterpart from section \ref{3}, and we consider the operator
$D^{\Cal W}:\Ga(H_0)\to\Ga(H_1)$. By conformal invariance, the system
$D^{\Cal W}(f)=0$ must be among the systems covered by Theorem
\ref{3.5}, and the above procedure identifies its solutions with
parallel sections of $\Cal W$. Since $\nabla^{\Cal W}$ is flat in the
locally conformally flat case, the space of parallel sections has
dimension $\dim(\Bbb W)$. Moreover, two solutions of the system
coincide if and only if their images under $L$ have the same value in
one point.

\section{Generalizations}\label{5}

In this last part, we briefly sketch how the developments of sections
\ref{3} and \ref{4} can be carried over to larger classes of geometric
structures. 

\subsection{The prolongation procedure for general $|1|$--graded Lie
  algebras}\label{5.1} 
The algebraic developments in \ref{3.1}--\ref{3.2a} generalize without
problems to a semisimple Lie algebra $\frak g$ endowed with a
$|1|$--grading, i.e.~a grading of the form $\frak
g=\fg_{-1}\oplus\fg_0\oplus\fg_1$. Given such a grading it is easy to
see that it is the eigenspace decomposition of $\ad(E)$ for a uniquely
determined element $E\in\frak g_0$. The Lie subalgebra $\fg_0$ is
automatically the direct sum of a semisimple part $\fg'_0$ and a
one--dimensional center spanned by $E$. This gives rise to
$E$--eigenspace decompositions for irreducible representations. Again
irreducible representations of $\fg$ may be parametrized by pairs
consisting of an irreducible representation of $\fg'_0$ and an integer
$\geq 1$. Then all the developments of \ref{3.1}--\ref{3.2a} work
without changes.

Next choose a Lie group $G$ with Lie algebra $\fg$ and let $G_0\subset
G$ be the subgroup consisting of those elements whose adjoint action
preserves the grading of $\fg$. Then this action defines an
infinitesimally effective homomorphism $G_0\to GL(\fg_{-1})$. In
particular, the semisimple part $G'_0$ of $G_0$ is a (covering of a)
subgroup of $GL(\fg_{-1})$, so this defines a type of geometric
structure on manifolds of dimension $\dim(\fg_{-1})$. This structure
is linked to representation theory of $G'_0$ in the same way as
Riemannian geometry is linked to representation theory of $O(n)$. 

For manifolds endowed with a structure of this type, there is an
analog of the prolongation procedure described in \ref{3.3}--\ref{3.5}
with closely parallel proofs, see \cite{BCEG}. The only change is that
instead of the Levi--Civita connection one uses any linear connection
on $TM$ which is compatible with the reduction of structure group.
There are some minor changes if this connection has torsion. The
systems that this procedure applies to are the following. One chooses
an irreducible representation $\Bbb W_0$ of $G'_0$ and an integer
$r\geq 1$. Denoting by $W_0$ the bundle corresponding to $\Bbb W_0$,
one then can handle systems whose principal symbol is (a multiple of)
the projection from $S^rTM\otimes W_0$ to the subbundle corresponding
to the irreducible component of maximal highest weight in
$S^r\fg_1\otimes\Bbb W_0$.

The simplest example of this situation is $\fg=\frak{sl}(n+1,\Bbb R)$
endowed with the grading $\begin{pmatrix}\fg_0 & \fg_1\\
  \fg_{-1}&\fg_0\end{pmatrix}$ with blocks of sizes $1$ and $n$. Then
$\fg_{-1}$ has dimension $n$ and $\fg_0\cong\frak{gl}(n,\Bbb R)$. For
the right choice of group, one obtains $G'_0=SL(n,\Bbb R)$, so the
structure is just a volume form on an $n$--manifold.

There is a complete description of $|1|$--gradings of semisimple Lie
algebras in terms of structure theory and hence a complete list of the
other geometries for which the procedure works. One of these is
related to almost quaternionic structures, the others can be described
in terms of identifications of the tangent bundle with a symmetric or
skew symmetric square of an auxiliary bundle or with a tensor product
of two auxiliary bundles.

\subsection{Invariant differential operators for
  AHS--structures}\label{5.2} 
For a group $G$ with Lie algebra
$\fg=\fg_{-1}\oplus\fg_0\oplus\fg_1$ as in \ref{5.1}, one defines
$P\subset G$ as the subgroup of those elements whose adjoint action
preserves the subalgebra $\fg_0\oplus\fg_1=:\frak p$. It turns out
that $\frak p$ is the Lie algebra of $P$ and $G_0\subset P$ can also
be naturally be viewed as a quotient of $P$.

On manifolds of dimension $\dim(\fg_{-1})$ we may consider reductions
of structure group to the group $G_0$. The passage from $G'_0$ as
discussed in \ref{5.1} to $G_0$ is like the passage from Riemannian to
conformal structures. As in \ref{4.2}, one may look at extensions of
the principal $G_0$--bundle defining the structure to a principal
$P$--bundle $\Cal G$ endowed with a normal Cartan connection
$\om\in\Om^1(\Cal G,\frak g)$. In the example $\fg=\frak{sl}(n+1,\Bbb
R)$ with $\fg_0=\frak{gl}(n,\Bbb R)$ from \ref{5.1}, the principal
$G_0$--bundle is the full frame bundle, so it contains no information.
One shows that such an extension is equivalent to the choice of a
projective equivalence class of torsion free connections on $TM$. In
all other cases (more precisely, one has to require that no simple
summand has this form) Cartan's result on conformal structures can be
generalized to show that such an extension is uniquely possible for
each given $G_0$--structure, see e.g. \cite{CSS2}.

\noindent
The structures equivalent to such Cartan connections are called
AHS--structures in the literature. Apart from conformal and projective
structures, they also contain almost quaternionic and almost
Grassmannian structures as well as some more exotic examples, see
\cite{CSS2}. For all these structures, the procedure from section
\ref{4} can be carried out without changes to construct differential
operators which are intrinsic to the geometry.

\subsection{More general geometries}\label{5.3}
The construction of invariant differential operators from section
\ref{4} applies to a much larger class of geometric structures. Let
$\frak g$ be a semisimple Lie algebra endowed with a $|k|$--grading,
i.e.~a grading of the form $\fg=\fg_{-k}\oplus\dots\oplus\fg_k$ for
some $k\geq 1$, such that $[\fg_i,\fg_j]\subset\fg_{i+j}$ and such
that the Lie subalgebra $\fg_-:=\fg_{-k}\oplus\dots\oplus\fg_{-1}$ is
generated by $\fg_{-1}$. For any such grading, the subalgebra $\frak
p:=\fg_0\oplus\dots\oplus\fg_k\subset\fg$ is a \textit{parabolic}
subalgebra in the sense of representation theory. Conversely, any
parabolic subalgebra in a semisimple Lie algebra gives rise to a
$|k|$--grading. Therefore, $|k|$--gradings are well understood and can
be completely classified in terms of the structure theory of semisimple
Lie algebras.

Given a Lie group $G$ with Lie algebra $\fg$ one always finds a closed
subgroup $P\subset G$ corresponding to the Lie algebra $\frak p$. The
homogeneous space $G/P$ is a so--called \textit{generalized flag
  variety}. Given a smooth manifold $M$ of the same dimension as
$G/P$, a \textit{parabolic geometry} of type $(G,P)$ on $M$ is given
by a principal $P$--bundle $p:\Cal G\to M$ and a Cartan connection
$\om\in\Om^1(\Cal G,\frak g)$.

In pioneering work culminating in \cite{Tanaka}, N.~Tanaka has shown
that assuming the conditions of regularity and normality on the
curvature of the Cartan connection, such a parabolic geometry is
equivalent to an underlying geometric structure. These underlying
structures are very diverse, but during the last years a uniform
description has been established, see the overview article
\cite{Srni05}. Examples of these underlying structures include
partially integrable almost CR structures of hypersurface type, path
geometries, as well as generic distributions of rank two in dimension
five, rank three in dimension six, and rank four in dimension seven.
For all these geometries, the problem of constructing differential
operators which are intrinsic to the structure is very difficult.

The BGG machinery developed in \cite{CSS:BGG} and \cite{CD}
offers a uniform approach for this construction, but compared to the
procedure of section \ref{4} some changes have to be made. One again
has a grading element $E$ which leads to an eigenspace decomposition
$\Bbb W=\Bbb W_0\oplus\dots\oplus\Bbb W_N$ of any finite dimensional
irreducible representation of $\fg$. As before, we have
$\fg_i\cdot\Bbb W_j\subset \Bbb W_{i+j}$. Correspondingly, this
decomposition is only invariant under a subgroup $G_0\subset P$ with
Lie algebra $\fg_0$, but each of the subspaces $\Bbb
W_i\oplus\dots\oplus\Bbb W_N$ is $P$--invariant. The theory of tractor
bundles and tractor connections works in this more general setting
without changes, see \cite{TAMS}.

Via the Cartan connection $\om$, the tangent bundle $TM$ can be
identified with $\Cal G\x_P(\frak g/\frak p)$ and therefore
$T^*M\cong\Cal G\x_P(\frak g/\frak p)^*$. Now the annihilator of
$\frak p$ under the Killing form of $\fg$ is the subalgebra $\frak
p_+:=\fg_1\oplus\dots\oplus\fg_k$. For a tractor bundle $\Cal W=\Cal
G\x_P\Bbb W$, the bundles of $\Cal W$--valued forms are therefore
associated to the representations $\La^k\frak p_+\otimes\Bbb W$. 

Since we are now working with the nilpotent Lie algebra $\frak p_+$
rather than with an Abelian one, we have to adapt the definition of
$\partial^*$. In order to obtain a differential, we have to add terms
which involve the Lie bracket on $\frak p_+$. The resulting map
$\partial^*$ is $P$--equivariant, and the quotients
$\ker(\partial^*)/\im(\partial^*)$ can be computed as representations
of $\fg_0$ using Kostant's theorem. As far as $\partial$ is concerned,
we have to identify $\frak g/\frak p$ with the nilpotent subalgebra
$\fg_-:=\fg_{-k}\oplus\dots\oplus\fg_{-1}$. Then we can add terms
involving the Lie bracket on $\fg_-$ to obtain a map $\partial$ which
is a differential. As the identification of $\fg/\frak p$ with
$\fg_-$, the map $\partial$ is not equivariant for the $P$--action but
only for the action of a subgroup $G_0$ of $P$ with Lie algebra $\fg_0$. 

The $P$--equivariant map $\partial^*$ again induces vector bundle
homomorphisms on the bundles of $\Cal W$--valued differential forms.
We can extend the normal tractor connection to the covariant exterior
derivative $d^\Cal W$. As in \ref{4.4}, the lowest homogeneous
component of $d^{\Cal W}$ is tensorial and induced by $\partial$ which
is all that is needed to get the procedure outlined in \ref{4.4}
going.  Also the results for structures which are locally isomorphic
to $G/P$ discussed in \ref{4.4} extend to general parabolic
geometries.

The question of analogs of the prolongation procedure from section
\ref{3} for arbitrary parabolic geometries has not been completely
answered yet. It is clear that some parts generalize without problems.
For other parts, some modifications will be necessary. In particular,
the presence of non--trivial filtrations of the tangent bundle makes
it necessary to use the concept of weighted order rather than the
usual concept of order of a differential operator and so on. Research
in this direction is in progress.


\begin{thebibliography}{10}
  
\bibitem{BEG} {\sc T.N. Bailey, M.G. Eastwood, A.R. Gover}, {\sl
    Thomas's structure bundle for conformal, projective and related
    structures}, Rocky Mountain J. \textbf{24} (1994), 1191--1217.
  
\bibitem{Baston} {\sc R.J. Baston}, {\sl Almost Hermitian symmetric
    manifolds, I: Local twistor theory; II: Differential invariants},
  Duke Math. J., \textbf{63} (1991), 81--111, 113--138.
  
\bibitem{BGG} {\sc I.N. Bernstein, I.M. Gelfand, S.I. Gelfand}, {\sl
    Differential operators on the base affine space and a study of
    $\mathfrak g$--modules}, in ``Lie Groups and their
  Representations'' (ed. I.M.  Gelfand) Adam Hilger 1975, 21--64.
  
\bibitem{BCEG} {\sc T. Branson, A. \v Cap, M.G. Eastwood , A.R. Gover},
  {\sl Prolongation of geometric overdetermined systems}, Internat. J. Math. 
\textbf{17} No. 6 (2006) 641--664, available online as math.DG/0402100 

\bibitem{CD}{\sc  D.M.J.Calderbank, T.Diemer}, {\sl Differential invariants and
  curved Bernstein-Gelfand-Gelfand sequences}, J. Reine Angew. Math.
  \textbf{537} (2001) 67--103.
  
\bibitem{Srni05} {\sc A. \v Cap}, {\sl Two constructions with parabolic
  geometries}, preprint math.DG/0504389, to appear in the proceedings
  of the 25th winter school ``Geometry and physics'', Suppl. Rend.
  Circ. Mat. Palermo.

\bibitem{TAMS} {\sc A. \v Cap, A.R. Gover}, {\sl Tractor Calculi for
Parabolic Geometries}, Trans. Amer. Math. Soc. \textbf{354} (2002),
1511-1548.
  
\bibitem{CSS2} {\sc A. \v Cap, J. Slov\'ak, V. Sou\v cek}, {\sl Invariant
  operators on manifolds with almost Hermitian symmetric structures,
  II. Normal Cartan connections}, Acta Math. Univ. Commenianae,
  \textbf{66} (1997), 203--220.

\bibitem{CSS:BGG} {\sc A. \v Cap, J. Slov\'ak, V. Sou\v cek},
{\sl Bernstein--Gelfand--Gelfand sequences}, Ann. of Math. \textbf{154}
no. 1 (2001) 97--113.

\bibitem{Mike:prolon}{\sc M.G. Eastwood}, {\sl Prolongations of linear
    overdetermined systems on affine and Riemannian manifolds},
  Rend. Circ. Mat. Palermo Suppl.  No. \textbf{75}  (2005), 89--108. 

\bibitem{Mike:symmetries}  {\sc M.G. Eastwood},  {\sl Higher symmetries of the
  Laplacian}, Ann. of Math. \textbf{161} no. 3 (2005) 1645--1665.
  
\bibitem{Eastwood-Rice} {\sc M.G. Eastwood, J.W. Rice}, {\sl
    Conformally invariant differential operators on Minkowski space
    and their curved analogues}, Commun. Math. Phys.  \textbf{109}
  (1987), 207--228.

\bibitem{Fegan} {\sc H.D. Fegan}, {\sl Conformally invariant first order
  differential operators}, Quart. J. Math. \textbf{27} (1976)
  371--378. 
  
\bibitem{Kostant} {\sc B. Kostant}, {\sl Lie algebra cohomology and the
  generalized Borel--Weil theorem}, Ann. of Math. \textbf{74} no. 2
  (1961), 329--387.
  
\bibitem{Lepowsky} {\sc J. Lepowsky},{\sl  A generalization of the
  Bernstein--Gelfand--Gelfand resolution}, J. of Algebra \textbf{49}
  (1977), 496--511.

\bibitem{Tanaka} {\sc N. Tanaka}, {\sl On the equivalence problem associated
  with simple graded Lie algebras}, Hokkaido Math. J., \textbf{8}
  (1979), 23--84.
  
\bibitem{Thomas} {\sc T.Y. Thomas}, {\sl On conformal geometry}, Proc. N.A.S.
  \textbf{12} (1926), 352--359; {\sl Conformal tensors}, Proc. N.A.S.
  \textbf{18} (1931), 103--189.


\end{thebibliography}
\end{document}